\documentclass[12pt,reqno]{amsart}

\usepackage{mathabx}
\usepackage{amssymb}
\usepackage{bbm}
\usepackage{mathrsfs}
\usepackage[all]{xy}
\usepackage{bm}
\usepackage{float}

\newtheorem{theo}{Theorem}[section]
\newtheorem{prop}[theo]{Proposition}
\newtheorem{lem}[theo]{Lemma}
\newtheorem{coro}[theo]{Corollary}

\theoremstyle{definition}
\newtheorem{defi}[theo]{Definition}

\theoremstyle{remark}
\newtheorem{rem}[theo]{Remark}
\newtheorem{rems}[theo]{Remarks}
\newtheorem{ex}[theo]{Example}
\newtheorem{conv}[theo]{Convention}

\newtheorem*{qu}{Question}

\newcommand{\br}{ }
\newcommand{\brr}{, }

\newcommand{\Gal}{\mathop{\text{\rm Gal}}\nolimits}
\newcommand{\Tr}{\mathop{\text{\rm Tr}}\nolimits}
\newcommand{\Sym}{\mathop{\text{\rm Sym}}\nolimits}
\newcommand{\AGL}{\mathop{\text{\rm AGL}}\nolimits}
\newcommand{\PGL}{\mathop{\text{\rm PGL}}\nolimits}
\newcommand{\res}{\mathop{\text{\rm res}}\nolimits}

\newcommand{\Br}{\mathop{\text{\rm Br}}\nolimits}
\newcommand{\Sp}{\mathop{\text{\rm Sp}}\nolimits}

\newcommand{\sep}{\text{\rm sep}}
\newcommand{\ds}{\text{\rm ds}}
\newcommand{\ev}{\text{\rm ev}}
\newcommand{\maxi}{\text{\rm max}}

\newcommand{\bbF}{{\mathbbm F}}
\newcommand{\bbQ}{{\mathbbm Q}}
\newcommand{\bbZ}{{\mathbbm Z}}

\newcommand{\calM}{{\mathscr{M}}}
\newcommand{\calO}{{\mathscr{O}}}
\newcommand{\calV}{{\mathscr{V}}}

\newcommand{\Ab}{{\text{\bf A}}}
\newcommand{\Pb}{{\text{\bf P}}}

\newcounter{abc}
\newenvironment{abc}{\begin{list}{\rm \alph{abc}) }%
{\usecounter{abc} \leftmargin=0.0pt \labelsep=0.0pt %
\listparindent=0.0pt \labelwidth=0.0pt \parsep=\smallskipamount%
 \itemsep=0.0pt \topsep=0.0pt \partopsep=\smallskipamount}}{\end{list}}

\newcounter{iii}
\newenvironment{iii}{\begin{list}{\rm \roman{iii}) }%
{\usecounter{iii} \leftmargin=0.0pt \labelsep=0.0pt %
\listparindent=0.0pt \labelwidth=0.0pt \parsep=\smallskipamount%
 \itemsep=0.0pt \topsep=0.0pt \partopsep=\smallskipamount}}{\end{list}}

\makeatletter
\def\hsmash{\relax 
  \ifmmode\def\next{\mathpalette\mathhsm@sh}\else\let\next\makehsm@sh
  \fi\next}
\def\makehsm@sh#1{\setbox\z@\hbox{#1}\finhsm@sh}
\def\mathhsm@sh#1#2{\setbox\z@\hbox{$\m@th#1{#2}$}\finhsm@sh}
\def\finhsm@sh{\wd\z@\z@ \box\z@}
\makeatother

\def\rightend#1#2{{%
 \leavevmode\nobreak\hskip .5em plus 1fil
 \penalty600 \hskip 0pt plus -1filll
 \vadjust{}\nobreak\hskip 0pt plus 1filll%
 #1\parfillskip=#2\relax \par}}

\def\eop{\ifmmode\rule[-22pt]{0pt}{1pt}\ifinner\tag*{$\square$}\else\eqno{\square}\fi\else\rightend{$\square$}{0pt}\fi}

\thanks{}

\title[On plane quartics with a Galois invariant Cayley octad]{On plane quartics with a Galois invariant \\Cayley octad}

\begin{document}

\author{Andreas-Stephan Elsenhans}

\address{Institut f\"ur Mathematik\\ Universit\"at W\"urzburg\\ Emil-Fischer-Stra\ss e 30\\ D-97074 W\"urzburg\\ Germany}
\email{stephan.elsenhans@mathematik.uni-wuerzburg.de}
\urladdr{https://math.uni-paderborn.de/ag/ca/elsenhans/}

\author{J\"org Jahnel}

\address{\mbox{Department Mathematik\\ \!Univ.\ \!Siegen\\ \!Walter-Flex-Str.\ \!3\\ \!D-57068 \!Siegen\\ \!Germany}}
\email{jahnel@mathematik.uni-siegen.de}
\urladdr{http://www.uni-math.gwdg.de/jahnel}

\thanks{{\em Acknowledgement.} We wish to thank the anonymous referee, whose remarks helped us to significantly improve this~article.}

\date{September~6,~2018.}



\begin{abstract}
We describe a construction of plane quartics with prescribed Galois operation on the 28~bitangents, in the particular case of a Galois invariant Cayley octad. As~an application, we solve the inverse Galois problem for degree two del Pezzo surfaces in the corresponding particular case.
\end{abstract}

\maketitle
\thispagestyle{empty}

\section{Introduction}

It is well known~\cite{Pl} that a nonsingular plane quartic
curve~$C$
over an algebraically closed field of
characteristic~$\neq \!2$
has exactly 28 bitangents. The~same is still true if the base field is only separably closed, as is easily deduced from~\cite[Theorem~1.6]{Va}.
If~$C$
is defined over a separably non-closed
field~$k$
then the bitangents are defined over a finite extension field
$l$
of~$k$,
which is normal and separable, and permuted by the Galois
group~$\Gal(l/k)$.

By~far not every permutation
in~$S_{28}$
may~occur.
It~is well-known (cf.,~e.g.,~\cite{PSV}) that the subgroup
$G \subset S_{28}$
of all admissible permutations is independent of the choice
of~$C$
and isomorphic to
$\Sp_6(\bbF_{\!2})$.
A~natural question arising is thus the~following. 

\begin{qu}
Given~a
field~$k$
and a subgroup
$g \subseteq G$,
does there exist a nonsingular plane
quartic~$C$
over~$k$,
for which the group homomorphism
$\Gal(k^\sep/k) \to G \subset S_{28}$,
given by the Galois operation on the 28 bitangents, has
image~$g$?
This~question depends only on the conjugacy class of the
subgroup~$g \subseteq G$.
\end{qu}

\begin{rem}
At~least when
$k$
is a global field, the Galois operation on the 28~bitangents and, in particular, the subgroup
of~$G$
being its image, form important invariants concerning the arithmetic
of the
curve~$C$.
For~interesting applications in the closely related situations of Del Pezzo surfaces of low degree, we advise the reader to consult Yu.\,I.~Manin's book~\cite{Ma} or the article \cite{Shi} of T.~Shioda.
\end{rem}

The group
$G \cong \Sp_6(\bbF_{\!2})$~is
the simple group of order
$1\,451\,520$.
It~has
$1369$~con\-ju\-gacy classes of subgroups. Among~these, there are eight maximal subgroups, which are of indices
$28$,
$36$,
$63$,
$120$,
$135$,
$315$,
$336$,
and
$960$,~respectively.

An~example of a nonsingular plane quartic
over~$\bbQ$
such that
$\Gal(\overline\bbQ/\bbQ)$
operates on the bitangents via the full
$\Sp_6(\bbF_{\!2})$
has been constructed by T.~Shioda~\cite[Section~3]{Shi} in 1993 and, almost at the same time, by R.~Ern\'e~\cite[Corollary~3]{Er}. Moreover,~there is an obvious approach to construct examples for the groups contained in the index
$28$~subgroup.
Indeed,~in this case, there is a rational~bitangent. One~may start with a cubic surface with the right Galois operation~\cite{EJ15}, blow-up a rational point, and use the connection between degree two del Pezzo surfaces and plane quartics~\cite[Theorem~3.3.5]{Ko96}, cf.\ the application discussed below (and mainly in Section~\ref{sec_drei}).
For~the group
$U_{63}$
of
index~$63$
and its subgroups, we gave a complete answer in our previous work~\cite{EJ17}.\smallskip

In~this article, we deal with the subgroup
$U_{36} \subset G$
of
index~$36$
and the groups contained~within. More precisely, we show the following result, which answers a more refined question than the one asked~above.

\begin{theo}
\label{general}
Let an infinite
field\/~$k$
of characteristic
not\/~$2$,
a normal and separable extension
field\/~$l$,
and an injective group~homomorphism
$$i\colon \Gal(l/k) \hookrightarrow U_{36}$$
be~given. Then there exists a nonsingular quartic
curve\/~$C$
over\/~$k$
such that\/
$l$
is the field of definition of the 28 bitangents and each\/
$\sigma \in \Gal(l/k)$
permutes the bitangents as described by\/
$i(\sigma) \in G \subset S_{28}$.
\end{theo}

Among the 1369 conjugacy classes of subgroups
of~$G \cong \Sp_6(\bbF_{\!2})$,
296 are contained
in~$U_{36}$.
However,~262 of these are also contained
in~$U_{63}$,
such that they are covered by~\cite{EJ17}, already, while 34 are~new.

The~subgroup
$U_{36}$
has a geometric meaning. The~assumption that the image of the natural homomorphism
$i\colon \Gal(k^\sep/k) \to G$
is contained
in~$U_{36}$
expresses the fact that there exists a Galois invariant Cayley octad (cf.~\cite[\S5]{Sen}, \cite[\S3]{PSV}, or \cite[Def\-i\-ni\-tion~6.3.4]{Do}) defining the
quartic~$C$.
Moreover,~one has an isomorphism of groups
$U_{36} \cong S_8$
and
$U_{36}$
permutes the eight points forming the octad~accordingly. We~use Cayley octads for our proof of existence, which is completely constructive.\smallskip

\begin{rems}
\label{mark}
\begin{iii}
\item
Theorem~\ref{general} is clearly not true, in general, when
$k$
is a finite~field. For~example, there cannot be a nonsingular quartic curve
over~$\bbF_{\!3}$,
all whose bitangents are
\mbox{$\bbF_{\!3}$-rational},
simply because the projective plane contains only 13 
$\bbF_{\!3}$-ra\-tio\-nal
lines.
\item
We~ignore about
characteristic~$2$
in this article, as this case happens to be very~different. Even~over an algebraically closed field, a plane quartic cannot have more than seven bitangents~\cite[p.~60]{SV}.
\end{iii}
\end{rems}

As~an application, one may answer the analogous question for degree two del Pezzo surfaces. The~double cover
of~$\Pb^2$,
ramified at a nonsingular quartic
curve~$C$
is a del Pezzo surface of degree two.
Here,~considerations can be made that are very similar to the ones~above. First~of all, it is well known~\cite[Theorem~26.2.(iii)]{Ma} that a del Pezzo
surface~$S$
of degree two over an algebraically closed field contains exactly 56 exceptional curves, i.e.\ such of self-intersection
number~$(-1)$.
Again,~the same is true when the base field is only separably~closed.
If~$S$
is defined over a separably non-closed
field~$k$
then the exceptional curves are defined over a normal and separable finite extension field
$l$
of~$k$
and permuted
by~$\Gal(l/k)$.
Once~again, not every permutation in
$S_{56}$
may~occur. The~maximal subgroup
$\smash{\widetilde{G} \subset S_{56}}$
that respects the intersection pairing is isomorphic to the Weyl
group~$W(E_7)$~\cite[Theorem~23.9]{Ma}.

Every~bitangent
of~$C$
is covered by exactly two of the exceptional curves
of~$S$.
Thus,~for the operation
of~$\smash{\Gal(k^\sep/k)}$
on the 56 exceptional curves
on~$S$,
there seem to be two independent conditions. On~one hand,
$\smash{\Gal(k^\sep/k)}$
must operate via a subgroup
of~$\smash{W(E_7) \cong \widetilde{G} \subset S_{56}}$.
On~the other hand, the induced operation on the blocks of size two must take place via a subgroup of
$\Sp_6(\bbF_{\!2}) \cong G \subset S_{28}$.
In~turns out, however, that there is an isomorphism
$\smash{W(E_7)/Z \stackrel{\cong}{\longrightarrow} \Sp_6(\bbF_{\!2})}$,
for
$Z \subset W(E_7)$
the centre, that makes the two conditions~equivalent. Cf.~\cite[Corollary~2.17]{EJ17}.

The~group
$\smash{\widetilde{G} \cong W(E_7)}$
already has
$8074$
conjugacy classes of~subgroups. Two~subgroups with the same image under the quotient~map
$\smash{p\colon \widetilde{G} \twoheadrightarrow \widetilde{G}/Z \stackrel{\cong}{\longrightarrow} G}$
correspond to del Pezzo surfaces of degree two that are quadratic twists of each~other. Theorem~\ref{general} therefore extends word-by-word to del Pezzo surfaces of degree two and homomorphisms
$\smash{\Gal(l/k) \hookrightarrow \widetilde{G}}$
with image contained in
$p^{-1}(U_{36})$.

\begin{rem}[Conjugacy classes in
$G$
versus conjugacy classes in
$U_{36}$]
As~noticed above, 296 of the conjugacy classes of subgroups of
$G \cong \Sp_6(\bbF_{\!2})$
are contained in
$S_8 \cong U_{36} \subset G$.
On~the other hand, one may readily check that
$S_8$
itself has precisely 296 conjugacy classes of~subgroups. In~other words, two subgroups
$U_1, U_2 \subseteq U_{36}$
that are conjugate
in~$G$
must be conjugate
in~$U_{36}$,~already.

This~result came to us originally as an experimental finding, and as quite a surprise. As~it is pure group theory and thus somewhat off the main topic of this article, we provide a proof in an~appendix.
\end{rem}

\begin{conv}
In~this article, by a {\em field,} we mean a field of
characteristic~$\neq\! 2$.
For~the convenience of the reader, the assumption on the characteristic will be repeated in the formulations of our final results, but not during the intermediate~steps.
\end{conv}

\noindent
{\em Computations.}
All computations are done using {\tt magma}~\cite{BCP}.

\section{Cayley octads with prescribed Galois operation}

All space quadrics over a
field~$k$
form a
$\Pb^9$,
the generic quadric being
$$
(T_0\, T_1\, T_2\, T_3)
\left(
\begin{array}{cccc}
a_{00} & a_{01} & a_{02} & a_{03} \\
a_{01} & a_{11} & a_{12} & a_{13} \\
a_{02} & a_{12} & a_{22} & a_{23} \\
a_{03} & a_{13} & a_{23} & a_{33}
\end{array}
\right)
\left(
\begin{array}{c}
T_0 \\ T_1 \\ T_2 \\ T_3
\end{array}
\right)
= 0 \, .$$
Among~these, the singular quadrics are parametrised by the quartic
$Q_s \subset \Pb^9$,
given by
$D_s = 0$,~for
$$D_s := \det
\left(
\begin{array}{cccc}
a_{00} & a_{01} & a_{02} & a_{03} \\
a_{01} & a_{11} & a_{12} & a_{13} \\
a_{02} & a_{12} & a_{22} & a_{23} \\
a_{03} & a_{13} & a_{23} & a_{33}
\end{array}
\right) .$$

\begin{lem}
\label{det_quart}
The quartic\/
$Q_s$
is singular, the singular points corresponding exactly to the quadrics of rank\/
$\leq\! 2$.\smallskip

\noindent
{\bf Proof.}
{\em
This~may be tested after base extension to the algebraic closure. Thus,~let us suppose that
$k = \overline{k}$.
The~quartic
$Q_s$
parametrises all singular quadrics, i.e.\ those of rank
$\leq\! 3$.
If~the rank of a
quadric~$q$
is exactly three then, after a linear transformation of coordinates, we may assume that
$q = T_0^2 + T_1^2 + T_2^2$.
At~the point where
$a_{00} = a_{11}= a_{22} = 1$
and all other coordinates vanish, the gradient
of~$D_s$
is~$a_{33}$.
Thus,~$q$
corresponds to a regular~point.

On the other hand, if the rank
of~$q$
is
$\leq\! 2$
then one may assume without restriction that
$q = T_0^2 + cT_1^2$,
for
$c = 1$
or~$0$.
At the point where
$a_{00} = 1$,
$a_{11} = c$,
and all other coordinates vanish, the gradient
of~$D_s$
is~zero. The~point is~singular.
}
\eop
\end{lem}

Let now
$q_1$,
$q_2$,
and~$q_3$
be three linearly independent quadratic forms, defined by the symmetric matrices
$M_1$,
$M_2$,
and~$M_3$.
These~span a~net
$$\Lambda := (Z(u_1q_1 + u_2q_2 + u_3q_3))_{(u_1:u_2:u_3) \in \Pb^2}$$
of quadrics, which contains its singular members at the~locus
$$C_\Lambda := Z(\det (u_1M_1 + u_2M_2 + u_3M_3)) \subset \Pb^2 \, .$$
This is a plane~quartic.

\begin{lem}
Assume that the base locus of the net\/
$\Lambda$
consists of eight distinct points, no four of which are~coplanar. Then~the quartic\/~$C_\Lambda$
is~nonsingular.\smallskip

\noindent
{\bf Proof.}
{\em
Again,~this may be tested after base change to the algebraic~closure. The~quartic
$C_\Lambda$
is then the intersection
of~$Q_s$
with a two-dimensional linear subspace
$E \subset \Pb^9$.
Thus,~there are two ways, in which
$C_\Lambda$
may become singular. Either,~the linear subspace intersects the singular locus
of~$Q_s$
or it meets
$Q_s$
somewhere~tangentially.

By~Lemma~\ref{det_quart}, the first option would mean that
$\Lambda = \langle q_1, q_2, q_3 \rangle$
contains a quadric of
rank~$\leq\! 2$.
Such~a quadric, however, splits into two linear forms, so the base locus
of~$\Lambda$,
given
in~$\Pb^3$
by
$q_1 = q_2 = q_3 = 0$,
is necessarily contained in the union of two~planes. This~contradicts our~assumptions.

Thus,~let us assume that
$E$
meets
$Q_s$
somewhere~tangentially. Without~restriction, the point of tangency corresponds to the space quadric
$q = T_0^2 + T_1^2 + T_2^2$.
The~tangent hyperplane at the corresponding point
on~$Q_s \subset \Pb^9$
is given
by~$a_{33} = 0$,
so that our assumption means that
$\Lambda = \langle q, q_1, q_2 \rangle$
and both
$q_1$
and~$q_2$
have coefficient zero
at~$T_3^2$.
But~then the base locus
of~$\Lambda$
contains
$(0\!:\!0\!:\!0\!:\!1)$
at least as a double point, which is a contradiction,~too.
}
\eop
\end{lem}

\begin{defi}
A set
$X \subset \Pb^3$
of eight points is called a {\em Cayley octad\/} if it is the base locus of a unique net of~quadrics. In~this case, we
write~$\Lambda_X$
for the net of quadrics
through~$X$.
\end{defi}

\begin{rems}
\begin{iii}
\item
Observe~that eight points in general position
in~$\Pb^3$
define only a pencil of quadrics. Thus,~Cayley octads are point sets in non-general position, although three general quadrics have eight points of~intersection.
\item
Let~$k$
be any field and
$X \subset \Pb^3_{k^\sep}$
be a Cayley octad that is invariant under
$\Gal(k^\sep/k)$.
Then~the net
$\Lambda_X$
is spanned by
\mbox{$k$-rational}
quadrics, so that
$C_{\Lambda_X}$
is a plane quartic defined
over~$k$.
\end{iii}
\end{rems}

\begin{prop}
\label{Cayley_bitangent}
Let\/~$k$
be any field and\/
$\smash{X \subset \Pb^3_{k^\sep}}$
be a Cayley octad that is invariant
under\/~$\Gal(k^\sep/k)$.
Assume~that the quartic\/
$C_{\Lambda_X}$
is nonsingular.\smallskip

\noindent
Then~the operation of\/
$\Gal(k^\sep/k)$
on the 28 bitangents
of\/~$C_{\Lambda_X}$
coincides with that on pairs of the eight points
of~$X$.
In~particular, if each of the eight points
of~$X$
is defined
over an extension
field\/~$l$
of\/~$k$
then the 28 bitangents
of\/~$C$
are defined
over\/~$l$.\smallskip

\noindent
{\bf Proof.}
{\em
As~is classically known (cf.~\cite[Proposition~5.4]{Sen} or \cite[Proposition~3.3]{PSV}), each pair of points
of~$X$
defines a bitangent
of~$C_{\Lambda_X}$,
which may be rationally computed from the~pair. The~assertion immediately follows from~this.
}
\eop
\end{prop}\pagebreak[3]

\begin{prop}
\label{Cayley_constr}
Let\/~$k$
be a field and\/
$f \in k[T]$
a separable polynomial of
degree\/~$8$,
whose coefficient
at\/~$T^7$
vanishes. Denote by\/
$l_0 := k[T]/(f)$
the \'etale algebra defined
by\/~$f$,
and let\/
$\alpha := (T \bmod f)$
be its natural~generator.\smallskip

\noindent
Furthermore,~write\/
$l$
for the splitting field
of\/~$f$,
and let\/
$\sigma_1,\ldots,\sigma_8\colon l_0 \rightarrow l$
be the eight homomorphisms, which are given by the roots
of\/~$f$.
The~operation
of\/~$\Gal(l/k)$
on the homomorphisms then yields an injection\/
$\iota\colon \Gal(l/k) \hookrightarrow S_8$.

\begin{abc}
\item
Then~the eight points\/
$(1\!:\!\sigma_i(\alpha)\!:\!\sigma_i(\alpha)^2\!:\!\sigma_i(\alpha)^4) \in \Pb^3(l)$,
for\/
$i=1,\ldots,8$,
form a Cayley
octad\/~$X$.
\item
Moreover,
$X$
is\/
\mbox{$\Gal(l/k)$-invariant},
the eight points being acted upon
by\/~$\Gal(l/k)$
as described
by\/~$\iota$.
\item
Suppose that, for each subset\/
$\{i_1, \ldots, i_4\} \subset \{1, \ldots, 8\}$
of size~four,
$$\sum_{j=1}^4 \sigma_{i_j}(\alpha) \neq 0 \, .$$
Then~the quartic\/
$C_{\Lambda_X}$
associated with\/
$X$
is~nonsingular.
\end{abc}\smallskip

\noindent
{\bf Proof.}
{\em
a)
The
\mbox{$2$-uple}
embedding maps
$(1\!:\!\sigma_i(\alpha)\!:\!\sigma_i(\alpha)^2\!:\!\sigma_i(\alpha)^4)$
to
\begin{equation}
\label{Veronese}
(1\!:\!\sigma_i(\alpha)^2\!:\!\sigma_i(\alpha)^4\!:\!\sigma_i(\alpha)^8\!:\!\sigma_i(\alpha)\!:\!\sigma_i(\alpha)^2\!:\!\sigma_i(\alpha)^4\!:\!\sigma_i(\alpha)^3\!:\!\sigma_i(\alpha)^5\!:\!\sigma_i(\alpha)^6) \, .
\end{equation}
We~first claim that the
$8\times10$-matrix
formed by the eight rows of the form (\ref{Veronese}) is of
rank~$7$.
For~this, we may reorder the columns, thereby omitting the repeated ones. The result is the
$8\times8$-matrix
$$
\left(
\begin{array}{cccccccc}
1 & \sigma_1(\alpha) & \sigma_1(\alpha)^2 & \sigma_1(\alpha)^3 & \sigma_1(\alpha)^4 & \sigma_1(\alpha)^5 & \sigma_1(\alpha)^6 & \sigma_1(\alpha)^8 \\
1 & \sigma_2(\alpha) & \sigma_2(\alpha)^2 & \sigma_2(\alpha)^3 & \sigma_2(\alpha)^4 & \sigma_2(\alpha)^5 & \sigma_2(\alpha)^6 & \sigma_2(\alpha)^8 \\
1 & \sigma_3(\alpha) & \sigma_3(\alpha)^2 & \sigma_3(\alpha)^3 & \sigma_3(\alpha)^4 & \sigma_3(\alpha)^5 & \sigma_3(\alpha)^6 & \sigma_3(\alpha)^8 \\
1 & \sigma_4(\alpha) & \sigma_4(\alpha)^2 & \sigma_4(\alpha)^3 & \sigma_4(\alpha)^4 & \sigma_4(\alpha)^5 & \sigma_4(\alpha)^6 & \sigma_4(\alpha)^8 \\
1 & \sigma_5(\alpha) & \sigma_5(\alpha)^2 & \sigma_5(\alpha)^3 & \sigma_5(\alpha)^4 & \sigma_5(\alpha)^5 & \sigma_5(\alpha)^6 & \sigma_5(\alpha)^8 \\
1 & \sigma_6(\alpha) & \sigma_6(\alpha)^2 & \sigma_6(\alpha)^3 & \sigma_6(\alpha)^4 & \sigma_6(\alpha)^5 & \sigma_6(\alpha)^6 & \sigma_6(\alpha)^8 \\
1 & \sigma_7(\alpha) & \sigma_7(\alpha)^2 & \sigma_7(\alpha)^3 & \sigma_7(\alpha)^4 & \sigma_7(\alpha)^5 & \sigma_7(\alpha)^6 & \sigma_7(\alpha)^8 \\
1 & \sigma_8(\alpha) & \sigma_8(\alpha)^2 & \sigma_8(\alpha)^3 & \sigma_8(\alpha)^4 & \sigma_8(\alpha)^5 & \sigma_8(\alpha)^6 & \sigma_8(\alpha)^8
\end{array}
\right) ,
$$
which is clearly of
rank~$\leq \!7$,
as the rightmost column is a linear combination of the other~seven. On~the other hand, the upper left
\mbox{$7\times7$-minor}
is Vandermonde and thus of~value
$$\prod_{1\leq i < j \leq 7} \!\![\sigma_i(\alpha) - \sigma_j(\alpha)] \neq 0 \, ,$$
which implies our~claim. Consequently,~the quadrics through
$X$
indeed form a
net~$\Lambda_X$.

In order to prove that
$X$
is indeed a Cayley octad, it now suffices to verify that the base locus
of~$\Lambda_X$
is zero dimensional. In~fact, in this case, according to Bezout, it cannot consist of more than the eight points~given.
To~show zero dimensionality, we first observe that the net
$\Lambda_X$
contains the
pencil~$L$
spanned by
$Z(q_1)$
and~$Z(q_2)$,
for~$q_1 := T_1^2-T_0T_2$
and~$q_2 := T_2^2-T_0T_3$.
The~pencil~$L$
contains the nonsingular~quadric
$Q := Z(q_1+q_2)$
and~has
$$C\colon T_1^2-T_0T_2 = T_2^2-T_0T_3 = 0$$
as its base locus, which is an irreducible~curve. Indeed,~on the affine chart
``$T_0 = 1$'',
one has the system of equations
$T_1^2=T_2$
and
$T_2^2=T_3$
that yields the rational parametrisation
$\Pb^1 \to C$,
$\smash{(s\!:\!t) \mapsto (s^4:s^3t:s^2t^2:t^4)}$,
for~$C$.
Consequently,~the base locus
of~$\Lambda_X$
could be of positive dimension only if the third generator
of~$\Lambda_X$
contained the
curve~$C$~entirely.

However,~left exactness of the global section functor, applied to the standard exact sequences
$$0 \longrightarrow \calO_{\Pb^3} \stackrel{\cdot (q_1+q_2)}{\longrightarrow} \calO_{\Pb^3}(2) \longrightarrow \calO_Q(2) \longrightarrow 0$$
and
$$0 \longrightarrow \calO_Q \stackrel{\cdot q_1}{\longrightarrow} \calO_Q(2) \longrightarrow \calO_C(2) \longrightarrow 0$$
shows that no quadratic form vanishes
on~$C$,
except for those linearly spanned
by~$q_1$
and~$q_2$.\smallskip

\noindent
b)
This~assertion is clear from the construction
of~$X$.\smallskip

\noindent
c)
Since~the eight points given are clearly distinct, all that is left to verify for 
$C_{\Lambda_X}$
being nonsingular is that no four of the points are~coplanar. But~this is clear,~too. Indeed, the~determinants
$$
\det
\left(
\begin{array}{cccc}
1 & \sigma_{i_1}(\alpha) & \sigma_{i_1}(\alpha)^2 & \sigma_{i_1}(\alpha)^4 \\
1 & \sigma_{i_2}(\alpha) & \sigma_{i_2}(\alpha)^2 & \sigma_{i_2}(\alpha)^4 \\
1 & \sigma_{i_3}(\alpha) & \sigma_{i_3}(\alpha)^2 & \sigma_{i_3}(\alpha)^4 \\
1 & \sigma_{i_4}(\alpha) & \sigma_{i_4}(\alpha)^2 & \sigma_{i_4}(\alpha)^4
\end{array}
\right) ,
$$
for
$1 \leq i_1 < \ldots < i_4 \leq 8$,
are easily calculated to~be
$$\big(\sigma_{i_1}(\alpha) + \sigma_{i_2}(\alpha) + \sigma_{i_3}(\alpha) + \sigma_{i_4}(\alpha)\big) \cdot \!\!\!\prod_{1\leq j_1 < j_2 \leq 4} \!\!\!\![\sigma_{j_1}(\alpha) - \sigma_{j_2}(\alpha)] \, ,$$
an expression, in which the factor on the left is nonzero by assumption while that on the right is nonzero by~construction.
}
\eop
\end{prop}

\begin{rem}
The~more obvious choice of the eight points
$(1\!:\!\sigma_i(\alpha)\!:\!\sigma_i(\alpha)^2\!:\!\sigma_i(\alpha)^3)$
does not lead to a Cayley~octad. There~is actually a net of quadrics through these eight points, spanned by
$Z(T_1^2-T_0T_2)$,
$Z(T_2^2-T_1T_3)$,
and~$Z(T_1T_2-T_0T_3)$,
but the base locus of this net is the twisted cubic~curve.
\end{rem}

We~may thus conclude the following~result.

\begin{theo}
\label{main_co}
Let an infinite
field\/~$k$
of characteristic
not\/~$2$,
a normal and separable extension
field\/~$l$,
and an injective group~homomorphism
$$i\colon \Gal(l/k) \hookrightarrow U_{36} \subset G \cong \Sp_6(\bbF_{\!2}) \subset S_{28}$$
be given. Then there exists a nonsingular quartic
curve\/~$C$
over\/~$k$
such that\/
$l$
is the field of definition of the 28 bitangents and each\/
$\sigma \in \Gal(l/k)$
permutes the bitangents as described by\/
$i(\sigma) \in G \subset S_{28}$.\smallskip

\noindent
{\bf Proof.}
{\em
The natural operation of
$S_8$
on the
$\smash{28 = \frac{8\cdot7}2}$
\mbox{2-sets} yields an injective group homomorphism
$S_8 \to S_{28}$,
the image of which is
exactly~$U_{36}$.
Thus,~the
homomorphism~$i$
given induces an injection
$\imath'\colon \Gal(l/k) \hookrightarrow S_8$.

This~immediately yields an \'etale
algebra~$l_0$
of degree eight. In~the case of a transitive subgroup,
$l_0$
is the field corresponding under the Galois correspondence to the point~stabiliser. In~general, one has to take the direct product of the fields, corresponding to the point stabilisers of the~orbits.

Moreover,~by the primitive element theorem, every element
of~$l_0$
is a generator, except for the union of finitely many
\mbox{$k$-linear}
subspaces of lower~dimension. The~same is still true for the trace zero subspace
$$h := \{\alpha \in l_0 \mid \Tr_{l_0/k}(\alpha) = 0\} \subset l_0 \, .$$
Indeed,~with
$\alpha \in l_0$,
the element
$\alpha' := \alpha - \frac18 \Tr_{l_0/k}(\alpha)$
of trace zero is a generator,~too.

Furthermore,~the eight homomorphisms
$\sigma_i\colon l_0\to l$
are
\mbox{$k$-linearly}
independent. Hence,~the subset
$$\textstyle h_c := \big\{\alpha \in h \mid \sum\limits_{j=1}^4 \!\sigma_{i_j}(\alpha) = 0 \text{ for some } 1 \leq i_1 < \ldots < i_4 \leq 8 \big\} \subset h \, ,$$
too, is a union of finitely many lower-dimensional
\mbox{$k$-linear}
subspaces.

Consequently,~as
$k$
is infinite, there exists some element
$\smash{\alpha \in h \setminus h_c}$
that is a generator of the \'etale
algebra~$l_0$.
We~take
$f \in k[T]$
to be the minimal
polynomial~of~$\alpha$.
Then,~one has
$l_0 \cong k[T]/(f)$
and the generator
$(T \bmod f)$
is mapped
to~$\alpha$
under this~isomorphism. The~coefficient
of~$f$
at~$T^7$
is zero, because of
$\smash{\Tr_{l_0/k}(\alpha) = 0}$.
Moreover,
$f$~is
separable, as, for
$\alpha$
a generator,
$\sigma_1(\alpha), \ldots, \sigma_8(\alpha)$
are automatically~distinct. Thus,~all the assumptions of Proposition~\ref{Cayley_constr} are~fulfilled.

Therefore,~we are given a Galois invariant Cayley octad
$X$
delivering a nonsingular plane
quartic~$C_{\Lambda_X}$.
Proposition~\ref{Cayley_bitangent} shows, finally, that the Galois operation on the bitangents
of~$C_{\Lambda_X}$
is exactly the one~desired.
}
\eop
\end{theo}

\begin{coro}[The case that
$k$
is a number field]

Let~$k$
be a number field and
$g$
a subgroup of
$G \cong \Sp_6(\bbF_{\!2})$
that is contained
in~$U_{36}$.
Then there exists a nonsingular quartic
curve~$C$
over~$k$
such that the natural permutation~representation
$$i\colon \Gal(\overline{k}/k) \longrightarrow G \subset S_{28}$$
on the 28 bitangents
of~$C$
has the
subgroup~$g$
as its~image.\smallskip

\noindent
{\bf Proof.}
{\em
According~to Theorem~\ref{main_co}, it suffices to show that, for every number
field~$k$
and each
subgroup~$g \subseteq U_{36}$,
there exists a normal extension
field~$l$
such that
$\Gal(l/k)$
is isomorphic
to~$g$.
This~is a particular instance of the inverse Galois problem for number fields and thus potentially~hard. But,~fortunately, the groups occurring are easy~enough.

In~fact, 296 of the 1369 conjugacy classes of subgroups
of~$G$
are contained
in~$U_{36}$.
These~are exactly the conjugacy classes of subgroups of
$U_{36} \cong S_8$.
It~does not happen that two subgroups
of~$U_{36}$
not being conjugate to each other
\mbox{become~conjugate~in}~$G$.
Cf.~Corollary~\ref{conj_U36}.

Among~these 296 conjugacy classes, 268 consist of solvable~groups. In~this case, the inverse Galois problem has been solved by I.\,R.~Shafarevich~\cite{Sha}, cf.~\cite[Theorem~9.5.1]{NSW}.
The~groups in the remaining 28 conjugacy classes turn out to be

\begin{iii}
\item
isomorphic to a symmetric group
$S_n$
or alternating group
$A_n$,
for~$n = 5$,
$6$,
$7$,
or~$8$,
in 13 cases,~and
\item
isomorphic to a direct product
$S_m \times S_n$,
$S_m \times A_n$,
or~$A_m \times A_n$,
in ten further~cases.
\item
The~remaining five cases are the simple group of
order~$168$,
occurring transitively and intransitively, the transitive subgroups of the types
$\PGL_2(\bbF_{\!7})$
and
$\AGL_3(\bbF_{\!2})$,
having the orders 
$336$
and~$1344$,
respectively, and the subdirect product of
$S_5$
and~$S_3$,
of
order~$360$.
\end{iii}

\noindent
The~groups of types i) and~ii) are classically known to occur as Galois groups over an arbitrary number
field~$k$
\cite[Theorem~9.4 and~9.5]{MM}. Moreover,~the~first three groups of type~iii) are realised
over~$\bbQ(t)$,
still
over~$\smash{\overline\bbQ(t)}$,
and consequently over
$k(t)$,
too, 
by the polynomials
$f_{7,3}$,
$f_{8,43}$,
and~$f_{8,48}$
from the tables in~\cite[Appendix]{MM}. Finally,~for the subdirect product, the~polynomial
$$[T^3 - 3(3t+1)T + 2(3t+1)] \cdot [T^5 - 5(-5t+1)T + 4(-5t+1)] \in \bbQ(t)[T]$$
enjoys the same~property. According~to the Hilbert irreducibility theorem~\cite[Section~9.2, Proposition~2]{Se}, in each case, there are infinitely many specialisations
of~$t$
to an element
of~$k$,
such that the resulting polynomial
over~$k$
has the same Galois~group.
}
\eop
\end{coro}

\begin{ex}
For~the polynomial
$$f(T) := T^8 + 42T^4 + 168T^2 + 1152T + 1197 \in \bbQ[T] \, ,$$
we find the nonsingular plane quartic
over~$\bbQ$,
given~by
\begin{align*}
2T_0^3T_2 - 5T_0^2T_2^2 - 3T_0T_1^3 - 6T_0T_1^2T_2 + 24T_0T_1T_2^2 - 8T_0T_2^3 + 6T_1^3T_2 + 12T_1^2T_2^2 & \\
{} + 24T_1T_2^3 - 16T_2^4 &= 0 \, .
\end{align*}
Here,~the Galois group 
of~$f$
is the simple group of
order~$168$,
realised as a transitive subgroup
of~$S_8$.
This~group is in fact doubly transitive, so that the operation on the 28~bitangents is transitive,~too.
\end{ex}

\begin{ex}
Over~the function field
$\bbF_{\!3}(t)$,
the~polynomial
\begin{align*}
f(T) := T^8 + 2tT^6 + 2t^2T^5 + (t^3 \!+\! 2t^2 \!+\! t + 2)T^4 + 2t^3T^3 + (2t^3 \!+\! t \!+\! 2)T^2 & \\[-1.5mm]
{} + (t^5 \!+\! t^4 \!+\! t^3 \!+\! 2t^2)T + (t^6 \!+\! t^4 \!+\! 2t^3 \!+\! t^2 \!+\! 1) & \in \bbF_{\!3}(t)[T]
\end{align*}
provides the same Galois~group. We~obtain the nonsingular plane quartic
over~$\bbF_{\!3}(t)$,
given~by
\begin{align*}
 & (t^{13} \!+\! t^{12} \!+\! 2t^{10} \!+\! t^8 \!+\! 2t^6 \!+\! 2t^4 \!+\! t^3 \!+\! t^2 \!+\! t \!+\! 1)T_0^4 \\[-1.5mm]
 & {} + (2t^9 \!+\! 2t^8 \!+\! t^7 \!+\! t^6 \!+\! 2t^5 \!+\! t^4 \!+\! 2t \!+\! 1)T_0^3T_1 \\[-1.5mm]
 & {} + (t^{10} \!+\! t^9 \!+\! t^8 \!+\! 2t^7 \!+\! t^6 \!+\! 2t^3 \!+\! t^2 \!+\! 2t \!+\! 1)T_0^3T_2
+ (t^4 \!+\! t^3 \!+\! 2t \!+\! 1)T_0^2T_1^2 \\[-1.5mm]
 & {} + (t^6 \!+\! t^5 \!+\! 2t^4 \!+\! 2t^3 \!+\! t \!+\! 1)T_0^2T_1T_2 + (t^7 \!+\! 2t^4 \!+\! t^3 \!+\! t^2 \!+\! 2t \!+\! 2)T_0^2T_2^2 + T_0T_1^3 \\[-1.5mm]
 & {} + tT_0T_1^2T_2 + (t^3 \!+\! 2t^2 \!+\! t \!+\! 2)T_0T_1T_2^2 + (2t^3 \!+\! t \!+\! 2)T_0T_2^3 + 2T_1T_2^3 = 0 \, .
\end{align*}
\end{ex}

\begin{ex}
For~the polynomial
$$f(T) := T^8 - 5T^6 - T^5 + 7T^4 + T^3 + 4T^2 + 1 \in \bbQ[T] \, ,$$
the nonsingular plane quartic
over~$\bbQ$,
given~by
\begin{align*}
T_0^4 - 2T_0^3T_1 - 5T_0^3T_2 + 6T_0^2T_1T_2 - 7T_0^2T_2^2 + 8T_0T_1^2T_2 - 6T_0T_1T_2^2 + 5T_0T_2^3 & \\
{} + 8T_1^3T_2 + T_1^2T_2^2 - 10T_1T_2^3 +& 2T_2^4 = 0
\end{align*}
results. Here,~the Galois group 
of~$f$
is~$A_8$.
The~operation on the 28~bitangents is transitive.
\end{ex}

\begin{ex}
Over~the function field
$\bbF_{\!3}(t)$,
the~polynomial
$$f(T) := T^8 + tT^5 + 2tT^2 + 1 \in \bbF_{\!3}(t)[T]$$
has Galois
group~$A_8$.
We~find the nonsingular plane quartic
over~$\bbF_{\!3}(t)$,
given~by
$$(t^2 + t)T_0^2T_1^2 + T_0T_1^3 + (t^2 + t)T_0^3T_2 + T_0^2T_1T_2 + 2tT_0T_2^3 + 2T_1T_2^3 = 0 \, . $$
\end{ex}

\begin{rem}
The examples above were chosen from the huge supply in the hope that they are of some particular~interest. For~instance, in both cases the groups occurring as Galois groups
are~simple. The~corresponding degree two del Pezzo surfaces are, independently of the twist chosen, all of Picard
rank~$1$.
\end{rem}

\begin{rem}[The situation of an algebraically closed base field---Moduli schemes]
There is a coarse moduli scheme of nonsingular quartic curves, provided by Geometric Invariant~Theory. It~is nothing but the quotient
$\calM := \calV/\PGL_3$,
for
$\calV \subset \Pb(\Sym^4(k^3)^*) \cong \Pb^{14}$
the open subscheme parametrising nonsingular plane quartics. Observe~here that, by~\cite[Remark~1.13]{Mu}, every nonsingular plane quartic corresponds to a
\mbox{$\PGL_3$-stable}
point
on~$\Pb(\Sym^4(k^3)^*)$.

The
\mbox{$\PGL_3$-invariants}
have been determined by J.~Dixmier~\cite{Di} and T.~Ohno, Ohno in his unpublished paper~\cite{Oh} proving~completeness. There~is a more recent treatment due to A.-S.~Elsenhans~\cite{El}, who also provided an implementation into~{\tt magma}. The~Dixmier--Ohno invariants yield an embedding
$$\calV/\PGL_3 \cong \calM \hookrightarrow \Pb(1,2,3,3,4,4,5,5,6,6,7,7,9)$$
into a weighted projective~space.

On~the other hand, there is the moduli
space~$\calM^\ev$
of projective equivalence classes of Cayley octads, cf.~\cite[Corollary~6.3.12]{Do}. As~is classically known, every nonsingular plane quartic may be obtained from a Cayley octad in 36 non-equivalent ways~\cite[\S14]{He}.
Accordingly,~$\calM^\ev$
naturally provides a
$36\!:\!1$
\'etale~covering
$$\pi\colon \calM^\ev \longrightarrow \calM$$
of the moduli
scheme~$\calM$
of nonsingular plane quartics~\cite[Theorem~6.3.19]{Do}.

Our~construction, from the geometric point of view, describes a particular kind of Cayley octads, namely those of the~type
$$\{ (1\!:\!x_i\!:\!x_i^2\!:\!x_i^4) \in \Pb^3(k) \mid i=1,\ldots,8\} \,,$$
for
$\{x_1, \ldots, x_8\} \subset k$
a subset consisting of eight elements such
that~$x_1 + \cdots + x_8 = 0$.
It~therefore yields a rational~map
$$\iota\colon \Ab^7 \dashrightarrow \calM^\ev \,.$$
Moreover,~a simple experiment in {\tt magma} reveals the fact that
$\iota$~is
dominant. Indeed, we verified dominance of the composition
$\pi \!\circ\! \iota\colon \Ab^7 \dashrightarrow \calM$
by showing that the tangent map
$T\pi_x\colon T_x \Ab^7 \to T_{\iota(x)} \calM$
at a single point
$x \in \Ab^7$
is~surjective. Here,~we used the first author's code~\cite{El} for calculating the Dixmier--Ohno invariants of a ternary quartic~form.

As~a consequence, we find that both,
$\calM^\ev$
and~$\calM$
are unirational~varieties. This~result, however, is not new. In~fact,
$\calM^\ev$
and~$\calM$
are even known to be~rational.
For~$\calM$,
this is due to P.~Katsylo~\cite[Theorem~0.1]{Ka}, while the rationality
of~$\calM^\ev$
follows from this in view of~\cite[Theorem~6.3.19]{Do}.

Anyway,~our approach to construct quartic curves with an arbitrarily given Galois operation on the 28~bitangents, having the image of Galois contained in
$U_{36} \subset G \cong \Sp_6(\bbF_{\!2})$,
seemed to run quite~smoothly. It~is our distinct impression that the rationality of the moduli
scheme~$\calM^\ev$
is the proper reason for~this.
\end{rem}

\section{Application to del Pezzo surfaces of degree two--Twisting}
\label{sec_drei}

There~is the double cover
$\smash{p\colon \widetilde{G} \to G}$
of finite groups, for
$\smash{W(E_7) \cong \widetilde{G} \subset S_{56}}$
and
$\Sp_6(\bbF_{\!2}) \cong G \subset S_{28}$,
which is given by the operation on the size two~blocks. The~kernel
of~$p$
is exactly the
centre~$\smash{Z \subset \widetilde{G}}$.
For~a subgroup
$\smash{H \subset \widetilde{G}}$,
one therefore has two~options.

\begin{iii}
\item
Either~$p|_H \colon H \to p(H)$
is two-to-one. Then
$H = p^{-1}(h)$,
for~$h := p(H)$.
In~this case,
$H$~contains
the centre
of~$\smash{\widetilde{G}}$
and, as abstract groups, one has an isomorphism
$H \cong p(H) \times \bbZ/2\bbZ$.
\item
Or~$p|_H \colon H \to p(H)$
is bijective.
\end{iii}
In~our geometric setting, the first case is the generic~one. More~precisely, let
$C\colon q=0$
be a nonsingular plane quartic such that the 28~bitangents are acted upon by the
group~$h \subseteq G$.
Then,~for
$\lambda$
an indeterminate, the 56~exceptional curves on
$S_\lambda\colon \lambda w^2 = q$
are operated upon
by~$h \times \bbZ/2\bbZ$.
For~particular choices
of~$\lambda$,
every subgroup
of~$\smash{\widetilde{G}}$
may be realised that has
image~$h$
under the
projection~$p$.

\begin{lem}
\label{twist}
Let~a
field\/~$k$
of characteristic
not\/~$2$,
a normal and separable extension
field\/~$l$,
and an injective group~homomorphism\/
$\smash{i\colon \Gal(l/k) \hookrightarrow \widetilde{G}}$
be~given.
Write\/~$l'$
for the subfield corresponding to\/
$i^{-1}(Z)$
under the Galois correspondence.

\begin{abc}
\item
Then there is a commutative diagram
$$
\xymatrix{
\Gal(l/k) \ar@{^{(}->}[rr]^{\;\;\;\;\;\;\;\;i} \ar@{^{}->>}[d]_{\res} && \widetilde{G} \ar@{^{}->>}[d]^p \\
\Gal(l'/k) \ar@{^{(}->}[rr]^{\;\;\;\;\;\;\;\;\overline\imath} && G \hsmash{ \, ,}
}
$$
the downward arrow on the left being the~restriction.
\item
Let\/~$C\colon q=0$
be a nonsingular plane quartic
over\/~$k$,
the 28 bitangents of which are defined
over\/~$l'$
and acted upon
by\/~$\Gal(l'/k)$
as described
by\/~$\overline\imath$.
Then~there exists some\/
$\lambda \in k^*$
such that the 56 exceptional curves of the degree two del Pezzo~surface
$$S_\lambda\colon \lambda w^2 = q$$
are defined
over\/~$l$
and each automorphism\/
$\sigma \in \Gal(l/k)$
permutes them as described by\/
$\smash{i(\sigma) \in \widetilde{G} \subset S_{56}}$.
\end{abc}\smallskip

\noindent
{\bf Proof.}
{\em
Cf.~\cite[Theorem~4.2]{EJ17}.
}
\eop
\end{lem}

\noindent
Thus, from our main result on plane quartics, we may draw the following conclusion.

\begin{theo}
\label{main_dP2}
Let an infinite
field\/~$k$
of characteristic
not\/~$2$,
a normal and separable extension
field\/~$l$,
and an injective group~homomorphism
$$i\colon \Gal(l/k) \hookrightarrow p^{-1}(U_{36})$$
be given. Then there exists a degree two del Pezzo
surface\/~$S$
over\/~$k$
such that\/
$l$
is the field of definition of the 56 exceptional curves and each\/
$\sigma \in \Gal(l/k)$
permutes them as described by\/
$\smash{i(\sigma) \in \widetilde{G} \subset S_{56}}$.\smallskip

\noindent
{\bf Proof.}
{\em
This follows from Lemma~\ref{twist} together with Theorem~\ref{main_co}.
}
\eop
\end{theo}

\begin{rem}[2-torsion Brauer classes]
There~is a unique~subgroup
$$H_{\Br} \subset p^{-1}(U_{36}) \subset \smash{\widetilde{G}} \subset S_{56}$$
of
index~$2$
that is intransitive of orbit
type~$[28,28]$.
Clearly,~one has
$H_{\Br} \cong U_{36} \cong S_8$.
The subgroup
$H_{\Br}$
is of
index~$72$
in~$\smash{\widetilde{G}}$.

Assume~that
$S$
is a degree two del Pezzo surface over a
field~$k$
such that the Galois group operating on the 56 exceptional curves is a subgroup
of~$H_{\Br}$.
In~this case,
$S$~carries
a global Brauer~class
$\alpha \in \Br(S_\lambda)_2$.
In~fact, the 2-torsion Brauer classes on degree two del Pezzo surfaces have been systematically studied by P.~Corn in~\cite{Co}. It~turns out that the conjugacy classes of subgroups
of~$W(E_7)$
that lead to such a Brauer class form a partially ordered set with exactly two maximal elements. The~subgroup
$H_{\Br}$
described above is one of them. It yields the Brauer classes of the second type in Corn's~terminology.

There~can be no doubt that this type is more interesting than the~other. A~short experiment shows the following facts.

\begin{iii}
\item
For 176 out of the 296 conjugacy classes
$H$
of subgroups
of~$H_{\Br}$,
the Brauer class
$\alpha$
is~nontrivial.
\item
For 87 of these, the Brauer class
$\alpha$
is non-cyclic. I.e.,~there is no normal subgroup
$H' \subset H$
with cyclic quotient annihilating the~class. Moreover,~somewhat surprisingly, in each of these 87 cases, the restriction to the
\mbox{$2$-Sylow}
subgroup
of~$H$
is still non-cyclic.

This~shows that, in the situation that
$k$
is a number field and
$\nu$
one of its places, the evaluation of such a Brauer class at a
\mbox{$k_\nu$-rational}
point cannot be computed in the usual, class field theoretic way (cf., e.g.,~\cite[\S6]{EJ10}).
\end{iii}
\end{rem}

\begin{rem}
Again, two subgroups
$U_1, U_2 \subseteq H_{\Br}$
not being conjugate to each other do not become conjugate when considered as subgroups
of~$\smash{\widetilde{G}}$.\smallskip

\noindent
{\bf Proof.}
Indeed,~observe the commutative diagram below,
$$
\xymatrix{
H_{\Br} \ar@{^{(}->}[rr] \ar@{^{}->}[d]_{p |_{H_{\Br}}}^{\cong} && \widetilde{G} \ar@{^{}->>}[d]^p \\
U_{36} \ar@{^{(}->}[rr] && G \hsmash{\,.}
}
$$
From this, one sees that if
$U_1, U_2 \subseteq H_{\Br}$
are conjugate
in~$\smash{\widetilde{G}}$
then
$p(U_1), p(U_2) \subseteq U_{36}$
are conjugate
in~$G$.
But then Corollary~\ref{conj_U36} shows that
$p(U_1), p(U_2) \subseteq U_{36}$
are already conjugate
in~$U_{36}$.
And,~finally, since
$p |_{H_{\Br}}$
is an isomorphism,
$U_1, U_2$
must be conjugate to each other
in~$H_{\Br}$.
\eop
\end{rem}

\section{Another application: Cubic surfaces with a Galois invariant double-six}

A~nonsingular cubic surface over an algebraically closed field contains exactly
27~lines. The~maximal subgroup
$\smash{G_\maxi \subset S_{27}}$
that respects the intersection pairing is isomorphic to the Weyl
group~$W(E_6)$~\cite[Theorem~23.9]{Ma}
of
order~$51\,840$.

A~double-six (cf.~\cite[Remark~V.4.9.1]{Ha} or \cite[Subsection~9.1.1]{Do}) is a configuration of twelve lines
$E_1, \ldots, E_6,E'_1, \ldots, E'_6$
such that

\begin{iii}
\item
$E_1,\ldots,E_6$
are mutually skew,
\item
$E'_1,\ldots,E'_6$
are mutually skew, and
\item
$E_i \!\cdot\! E'_j = 1$,
for\/
$i \neq j$,
$1 \leq i,j \leq 6$,
and\/
$E_i \!\cdot\! E'_i = 0$
for~$i = 1, \ldots, 6$.
\end{iii}

\noindent
Every cubic surface contains exactly 36 double-sixes, which are transitively acted upon
by~$G_\maxi$
\cite[Theorem~9.1.3]{Do}. Thus,~there is an index-36 subgroup
$U_\ds \subset G_\maxi$
stabilising a double-six. This~is one of the maximal subgroups
of~$G_\maxi \cong W(E_6)$.
Up~to conjugation,
$W(E_6)$
has maximal subgroups of indices
$2$,
$27$,
$36$,
$40$,
$40$,
and~$45$.

As~an application of Theorem~\ref{main_co} and Lemma~\ref{twist}, we now have a simpler proof for the following result, which was shown in~\cite[Theorem~5.1]{EJ17}.

\begin{theo}
\label{dsix}
Let an infinite
field\/~$k$
of characteristic
not\/~$2$,
a normal and separable extension
field\/~$l$,
and an injective group~homomorphism
$$i\colon \Gal(l/k) \hookrightarrow U_\ds \subset G_\maxi \cong W(E_6) \subset S_{27}$$
be~given. Then~there exists a nonsingular cubic
surface\/~$S$
over\/~$k$
such~that

\begin{iii}
\item
the 27 lines
on\/~$S$
are defined
over\/~$l$
and each\/
$\sigma \in \Gal(l/k)$
permutes them as described by\/
$i(\sigma) \in U_\ds \subset S_{27}$.
\item
$S$
is\/
\mbox{$k$-unirational}.
\end{iii}\smallskip

\noindent
{\bf Proof.}
{\em
There~is an injective homomorphism
$\iota\colon W(E_6) \hookrightarrow W(E_7)$
that corresponds to the blow-up of a point on a cubic~surface. The~subgroup
$\iota(U_\ds) \subset W(E_7)$
is of index
$36 \!\cdot\! 56 = 2016$.
It~is sufficient to show that the image
$\overline\iota(U_\ds)$
under the~composition
$$\overline\iota\colon W(E_6) \stackrel\iota\hookrightarrow W(E_7) \stackrel{p}\twoheadrightarrow W(E_7)/Z \cong G$$
is contained
in~$U_{36}$.
Indeed,~then Theorem~\ref{main_dP2} yields a degree two del Pezzo surface
$S'$
of degree two with a
\mbox{$k$-rational}
line~$L$
and the proposed Galois operation on the 27 lines that do not
meet~$L$.
Blowing~down
$L$,
we obtain a nonsingular cubic
surface~$S$
satisfying~i).
As~there is a
\mbox{$k$-rational}
point
on~$S$,
the blow down
of~$L$,
$S$~is
\mbox{$k$-unirational}~\cite[Theorem~2]{Ko02}.

To~show the group-theoretic claim, we first observe that the homomorphism
$\overline\iota$
is still an~injection. Furthermore,~the subgroup
$\overline\iota(W(E_6)) \subset G$
of
index~$28$
is exactly the point stabiliser of the transitive group
$G \subset S_{28}$.
The~subgroup
$U_{36} \subset S_{28}$
is still transitive, as
$U_{36} \cong S_8$,
operating on 2-sets, and
$S_8$
is doubly~transitive.
Thus,~$U_{36} \cap \overline\iota(W(E_6))$
is of
index~$28$
in~$U_{36}$.
Moreover,~$\overline\iota$
yields an isomorphism
$$\overline\iota |_{\overline\iota^{-1}(U_{36})}\colon \overline\iota^{-1}(U_{36}) \stackrel{\cong}{\longrightarrow} U_{36} \cap \overline\iota(W(E_6)) \, .$$
In~particular,
$\overline\iota^{-1}(U_{36})$
is of the same order as
$U_{36} \cap \overline\iota(W(E_6))$
and hence of
index~$36$
in the
group~$W(E_6)$.
Therefore,~up to conjugation,
$\overline\iota^{-1}(U_{36}) = U_\ds$,
which implies
$\overline\iota(U_\ds) = \overline\iota(\overline\iota^{-1}(U_{36})) \subseteq U_{36}$,
as~required.
}
\eop
\end{theo}

\appendix
\section{A group-theoretic observation}

\begin{prop}
\label{conj_subgr}
Let\/
$G$
be a finite group and\/
$H \subset G$
a subgroup. Consider~the permutation representation
\begin{eqnarray*}
\iota\colon G & \longrightarrow & \Sym(G/H) \,, \\
        g &\mapsto          & (aH \mapsto gaH)
\end{eqnarray*}
on left cosets
modulo\/~$H$,
and assume that, for all\/
$\omega, \omega' \in G/H$,
$\omega' \neq \omega$,
the following two conditions are~satisfied.

\begin{iii}
\item
There is some\/
$g \in G$
such that\/
$\iota(g) \in \Sym(G/H)$
interchanges\/
$\omega$
and\/
$\omega'$.
I.e.,~such that the permutation
$\iota(g)$
contains the
2-cycle~$(\omega \omega')$
in its canonical decomposition into disjoint~cycles.
\item
The twofold point stabiliser\/
$(G_\omega)_{\omega'}$
is a group having only inner automorphisms.
\end{iii}

\noindent
Then the following is~true. If two subgroups\/
$U_1, U_2 \subseteq H$
are conjugate
in\/~$G$
then they are conjugate already
in\/~$H$.\smallskip

\noindent
{\bf Proof.}
{\em
Consider, at first, two elements
$\omega, \omega' \in G/H$
such that
$\omega' \neq \omega$.
Then,~according to assumption~i), there is some
$g \in G$
such that
$\iota(g) = (\omega \omega')\cdots\;$.
Conjugation by
$g$
then clearly provides an automorphism
of~$(G_\omega)_{\omega'}$.
Moreover,~by assumption~ii), that must be inner. I.e.,~there exists some
$h \in (G_\omega)_{\omega'}$
satisfying
$hxh^{-1} = gxg^{-1}$
for every
$x \in (G_\omega)_{\omega'}$,
and therefore
$xh^{-1}g = h^{-1}gx$.
In other words, we found an~element
$$c_{\omega, \omega'} := h^{-1}g \in C_G\big((G_\omega)_{\omega'}\big)$$
in the centraliser
of~$(G_\omega)_{\omega'}$
that is non-trivial, because~of
$i(c_{\omega, \omega'}) = (\omega \omega')\cdots\;$.

Now suppose that
\begin{equation}
\label{conj}
u U_1 u^{-1} = U_2
\end{equation}
for some
$u \in G$.
We~let
$\omega := H$,
i.e.~the trivial~coset. Because~of
$U_1, U_2 \subseteq H$,
$\omega$~is
a fixed point of
$U_1$,
as well as
of~$U_2$.
Moreover,~assumption~(\ref{conj}) shows that
$\omega' := u^{-1}(\omega)$
must be a fixed point of
$U_1$,~too.

If~$\omega' = \omega$
then
$u \in H$
and the assertion is trivially~true. Otherwise,
$U_1 \subseteq (G_\omega)_{\omega'}$
and one has
\begin{eqnarray*}
U_2 & = & u U_1 u^{-1} \\
    & = & u c_{\omega, \omega'} U_1 c_{\omega, \omega'}^{-1} u^{-1} \,,
\end{eqnarray*}
since
$c_{\omega, \omega'}$
is a central~element. Here,
$u c_{\omega, \omega'}$
fixes
$\omega$,
which shows that
$u c_{\omega, \omega'} \in H$.%
}%
\qed
\end{prop}

\begin{rems}
\begin{iii}
\item
The~permutation representation
$\iota$
is always~transitive. Thus,~in the assumptions~i) and~ii), it suffices to restrict considerations to the case that
$\omega := H$.
\item
It is, of course, possible to find a purely group-theoretic formulation, avoiding~permutations. For~instance, the twofold point stabiliser of
$\omega = H$
and
$\omega' = xH$,
for~$x \not\in H$,
is~$H \cap xHx^{-1}$.
\end{iii}
\end{rems}

\begin{rem}
A~relevant particular case occurs when
$\iota$
doubly~transitive. Indeed,~in this case, assumption~i) is automatically fulfilled. Furthermore,~all twofold point stabilisers\/
$(G_\omega)_{\omega'}$
are isomorphic to each~other.
\end{rem}

\begin{coro}
\label{conj_U36}
If two subgroups\/
$U_1, U_2 \subseteq U_{36}$
are conjugate
in\/~$G \cong \Sp_6(\bbF_{\!2})$
then they are conjugate already
in\/~$U_{36}$.\smallskip

\noindent
{\bf Proof.}
{\em
In our case,
$G \cong \Sp_6(\bbF_{\!2})$
and 
$H = U_{36} \cong S_8$.
As~this is a maximal subgroup, the permutation representation
$\iota$
is primitive of
degree~$36$.
It~is well-known that
$\Sp_6(\bbF_{\!2})$
allows only one primitive permutation representation in this degree, and that this is doubly transitive~\cite[Theorem 7.7.A, for~$n=3$]{DM}.

Moreover,~$G_\omega \cong S_8$
allows only one conjugacy class of subgroups of
index~$35$,
and that contains the imprimitive wreath product
$S_4 \wr S_2$.
Thus,~all the twofold point stabilisers
of~$\iota$
are isomorphic
to~$S_4 \wr S_2$.
Finally,~this group allows only one faithful permutation representation in
degree~$8$
and, hence, has only inner automorphisms.
Thus,~Proposition~\ref{conj_subgr} applies and yields the~assertion.
}
\eop
\end{coro}

\setlength\parindent{0mm}
\end{document}